\documentclass[12pt,a4paper]{article}
\usepackage[T1]{fontenc}
\usepackage{mathptmx}
\usepackage{courier}
\usepackage[dvips]{graphicx}
\usepackage{psfrag}
\usepackage{amsmath,amssymb}        
\usepackage{latexsym}
\usepackage{url}

\def\D{\ensuremath{D}}
\def\www{.9\textwidth}

\title{A Non-Iterative Transformation Method for\\ Blasius Equation with Moving Wall or\\
 Surface Gasification}
\author{Riccardo Fazio \\
Department of Mathematics and Computer Science \\
University of Messina \\
Viale F. Stagno D'Alcontres, 31 \\
98166 Messina, Italy \\
E-mail: rfazio@unime.it \\
Home-page: http://mat521.unime.it/fazio}
\date{\today}
\begin{document}
\maketitle
\begin{abstract}
We define a non-iterative transformation method for Blasius equation with moving wall or surface gasification.
The defined method allows us to deal with classes of problems in boundary layer theory that, depending on a parameter, admit multiple or no solutions.
This approach is particularly convenient when the main interest is on the behaviour of the considered models with respect to the involved parameter.
The obtained numerical results are found to be in good agreement with those available in literature.    
\end{abstract}
\bigskip

\noindent
{\bf Key Words.} 
Blasius equation, non-iterative transformation method, BVPs on infinite intervals, moving wall, surface gasification. 
\bigskip

\noindent
{\bf AMS Subject Classifications.} 65L10, 34B40, 76D10, 76M55.


\section{Introduction.}
The problem of determining the steady two-dimensional motion of a
fluid past a flat plate placed edge-ways to the stream was
formulated in general terms, according to the boundary layer theory,
by Prandtl \cite{Prandtl:1904:UFK}, and was investigated in detail by
Blasius \cite{Blasius:1908:GFK}. 
The engineering interest was to calculate the shear
at the plate (skin friction),
which leads to the determination of the viscous drag on the plate, see for instance Schlichting 
\cite{Schlichting:2000:BLT}).

The celebrated Blasius problem is given by
\begin{align}\label{eq:Blasius} 
& {\displaystyle \frac{d^3 f}{d \eta^3}} + P \; f
{\displaystyle \frac{d^{2}f}{d\eta^2}} = 0 \nonumber \\[-1ex]
&\\[-1ex]
& f(0) = {\displaystyle \frac{df}{d\eta}}(0) = 0, \qquad
{\displaystyle \frac{df}{d\eta}}(\eta) \rightarrow 1 \quad \mbox{as}
\quad \eta \rightarrow \infty \ , \nonumber 
\end{align}
where $f$ and $\eta$ are suitable similarity variables and in literature we can find either $P = 1/2$ or $P=1$.
This is a boundary value problem (BVP) defined on the semi-infinite interval $[0, \infty)$.
According to Weyl \cite{Weyl:1942:DES}, the unique solution of (\ref{eq:Blasius}) has a positive second order derivative, which is monotone decreasing on $ [0, \infty) $ and approaches to zero as $ \eta $ goes to infinity.
The governing differential equation and the two boundary conditions
at $ \eta = 0$ in (\ref{eq:Blasius}) are invariant with respect to the scaling group of transformations
\begin{eqnarray}
& \eta^* = \lambda^{-\alpha}, \qquad f^* = \lambda^{\alpha} f
\label{eq:scaling:Blasius}
\end{eqnarray}
where $ \alpha $ is a nonzero constant: T\"opfer used $ \alpha = 1/3
$, see \cite{Topfer:1912:BAB}, but we have always put $ \alpha = 1 $ in order to simplify the
analysis. 
The mentioned invariance property has both analytical
and numerical interest. 
From a numerical viewpoint a non-iterative transformation method (ITM)
reducing the solution of (\ref{eq:Blasius}) to the solution of a related initial value problem (IVP) was defined by T\"opfer \cite{Topfer:1912:BAB}. 
Owing to that transformation, a simple existence and uniqueness Theorem
was given by J. Serrin \cite{Serrin:1970:ETS} as reported by Meyer \cite[pp. 104-105]{Meyer:1971:IMF} or Hastings and McLeod \cite[pp. 151-153]{Hastings:2012:CMO}. 
Let us note here that the mentioned invariance property is essential to the error analysis of the truncated boundary solution due to Rubel \cite{Rubel:1955:EET}, see Fazio \cite{Fazio:2002:SFB}.

Our main interest here is to extend T\"opfer's method to classes of problems in boundary layer theory involving a physical parameter.
This kind of extension was considered first by Na \cite{Na:1970:IVM}, see also Na \cite[Chapters 8-9]{Na:1979:CME}.
The application of a non-ITM to the Blasius equation with slip boundary condition, arising within the study of gas and liquid flows at the micro-scale regime \cite{Gad-el-Hak:1999:FMM,Martin:2001:BBL}, was considered already in \cite{Fazio:2009:NTM}.
Here we define a non-ITM for Blasius equation with moving wall considered by Ishak et al. \cite{Ishak:2007:BLM} or surface gasification studied by Emmons \cite{Emmons:1956:FCL} and recently by Lu and Law \cite{Lu:2014:ISB}. 
In particular, we find a way to solve non-iteratively the Sakiadis problem
\cite{Sakiadis:1961:BLBa,Sakiadis:1961:BLBb}.
For the solution of the Sakiadis problem by an ITM see Fazio \cite{Fazio:2015:ITM}.
The defined method allows us to deal with classes of problems in boundary layer theory that, depending on a parameter, admit multiple or no solutions.
This approach is particularly convenient when the main interest is on the behaviour of the considered models with respect to the involved parameter.    

\section{Moving wall}
According to Ishak et al. \cite{Ishak:2007:BLM} the differential problem governing a moving wall, with suitable boundary conditions, is given by
\begin{align}\label{eq:MWall}
& {\displaystyle \frac{d^3 f}{d \eta^3}} + \frac{1}{2} \; f
{\displaystyle \frac{d^{2}f}{d\eta^2}} = 0 \nonumber \\[-1ex]
&\\[-1ex]
&f(0) = 0 \ , \qquad
\frac{df}{d\eta} (0) = P  \ , \qquad
{\displaystyle \frac{df}{d\eta}}(\eta) \rightarrow 1-P \quad \mbox{as}
\quad \eta \rightarrow \infty \ , \nonumber
\end{align}
where $ P $ is a non-dimensional parameter given by the ratio of the wall to the flow velocities. 
Blasius problem (\ref{eq:Blasius}) is recovered from (\ref{eq:MWall}) by setting $P=0$.

\subsection{The non-ITM} 
The applicability of a non-ITM to the Blasius problem (\ref{eq:Blasius}) is a consequence of both: the invariance of the governing differential equation and the two boundary conditions at $\eta = 0$, and the non invariance of the asymptotic boundary condition under the scaling transformation (\ref{eq:scaling:Blasius}).
In order to apply a non-ITM to the BVP (\ref{eq:MWall}) we consider $ P $ as a parameter involved in the scaling invariance, i.e., we define the extended scaling group
\begin{equation}\label{eq:scalinv:wall}
f^* = \lambda f \ , \qquad \eta^* = \lambda^{-1} \eta \ , \qquad 
P^* = \lambda^{2} P \ .   
\end{equation}
Let us notice that, due to the given second boundary condition at $\eta = 0$ and the asymptotic boundary condition in (\ref{eq:MWall}), $P$ has to be transformed under the scaling group (\ref{eq:scalinv:wall}) with the same law of $\frac{df}{d\eta}(\eta)$.
By setting a value of $P^*$, we can integrate the Blasius equation in (\ref{eq:MWall}) written in the star variables on $[0, \eta^*_\infty]$, where $\eta^*_\infty$ is a suitable truncated boundary, with initial conditions
\begin{equation}\label{eq:ICs2}
f^*(0) = 0 \ , \quad \frac{df^*}{d\eta^*}(0) = P^* \ , \quad \frac{d^2f^*}{d\eta^{*2}}(0) = \pm 1 \ ,
\end{equation}
in order to compute an approximation $\frac{df^*}{d\eta^*}(\eta^*_\infty)$ for $\frac{df^*}{d\eta^*}(\infty)$ and the corresponding value of $\lambda$ according to the equation 
\begin{equation}\label{eq:lambda:wall}
\lambda = \left[ \frac{d f^*}{d \eta^{*}}(\eta^*_\infty)+P^* \right]^{1/2} \ .   
\end{equation} 
Once the value of $\lambda$ has been computed, by equation (\ref{eq:lambda:wall}), we can find the missed initial conditions 
\begin{equation}\label{eq:MIC}
\frac{df}{d\eta}(0) = \lambda^{-2} P^* \ , \qquad \frac{d^2f}{d\eta^{2}}(0) =  \lambda^{-3}\frac{d^2f^*}{d\eta^{*2}}(0) \ .
\end{equation}
Moreover, the numerical solution of the original BVP (\ref{eq:MWall}) can be computed by rescaling the solution of the IVP.
In this way we get the solution of a given BVP by solving a related IVP.

We remark here that the plus (for $P<0.5$) or minus (when $P>0.5$) sign must be used for the second derivative in (\ref{eq:ICs2}).
Moreover, the computation of a value at infinity is unsuitable from a numerical viewpoint and therefore we use a truncated boundary $\eta^*_{\infty}$ instead of infinity.
For the application of the method defined above, depending on the behaviour of the numerical solution, we used $\eta^*_{\infty} = 10$ or $\eta^*_{\infty} = 15$. 

In table \ref{tab:PMwall} we list sample numerical results obtained by the non-ITM for several values of $P^*$.
\begin{table}[!htb]
\renewcommand\arraystretch{1.1}
	\centering
		\begin{tabular}{cr@{.}lr@{.}lr@{.}lr@{.}l}
\hline \\[-3ex]
{${\displaystyle \frac{d^2f^*}{d{\eta^*}^2}(0)}$}
& \multicolumn{2}{c}%
{$ P^* $}
& \multicolumn{2}{c}%
{${\displaystyle \frac{df^*}{d\eta^*}(\infty)}$}
& \multicolumn{2}{c}%
{${\displaystyle \frac{d^2f}{d\eta^2}(0)}$}
& \multicolumn{2}{c}%
{$ P $} \\[2ex]
\hline
$1$ &
$-500$ & & $1$ & $55\D{04}$ & $5$ & $46\D{-07}$ & $-0$ & $033393$ \\
&$-100$ & & $2$ & $34\D{03}$ & $9$ & $42\D{-06}$ & $-0$ & $044591$ \\
&$-5$ &     & $36$ & $325698$ & $0$ & $005704$ & $-0$ & $159613$ \\
&$-1$ & $5$  & $4$ & $368544$ & $0$ & $205830$ & $-0$ & $522913$ \\
&$-1$ & $25$ & $3$ & $529165$ & $0$ & $290627$ & $-0$ & $548447$ \\
&$-1$ &      & $2$ & $917762$ & $0$ & $376537$ & $-0$ & $521441$ \\
&$-0$ & $75$ & $2$ & $503099$ & $0$ & $430814$ & $-0$ & $427814$ \\
&$-0$ & $5$  & $2$ & $250439$ & $0$ & $431797$ & $-0$ & $285643$ \\
&$0$ & & $2$ & $085393$ & $0$ & $332061$ & $0$ &  \\
&$1$ & & $2$ & $440648$ & $0$ & $156689$ & $0$ & $290643$ \\
&$5$ & & $5$ & $771518$ & $0$ & $028287$ & $0$ & $464187$ \\
&$100$ & & $1$ & $00\D{02}$ & $3$ & $53\D{-04}$ & $0$ & $499557$ \\
&$500$ & & $5$ & $00\D{02}$ & $3$ & $16\D{-05}$ & $0$ & $499960$ \\
$-1$ &
$100$ & & $99$ & $822681$ & $-3$ & $54\D{-04}$ & $0$ & $500444$ \\
&$10$ & & $9$ & $433763$ & $-0$ & $011673$ & $0$ & $514568$ \\
&$5$  & & $4$ & $182424$ & $-0$ & $035939$ & $0$ & $544519$ \\
&$2$  & & $0$ & $528464$ & $-0$ & $248722$ & $0$ & $790994$ \\
&$1$  & $719$ & $-4$ & $73\D{-05}$ & $-0$ & $443715$ & $1$ & $000027$ \\
\hline			
		\end{tabular}
	\caption{Moving wall boundary condition: non-ITM numerical results.}
	\label{tab:PMwall}
\end{table}
Here the $\D$ notation indicates that these results were computed in double precision.
As mentioned before, the case $P^*=P=0$ is the Blasius problem (\ref{eq:Blasius}).
In this case our non-ITM becomes the original method defined by T\"opfer \cite{Topfer:1912:BAB}.
For the Blasius problem, the obtained skin friction coefficient is in good agreement with the values available in literature, see for instance the value $0.332057336215$ computed by Fazio \cite{Fazio:1992:BPF} or the value $0.33205733621519630$, believed to be correct to all the sixteen decimal places, reported by Boyd \cite{Boyd:1999:BFC}. 
The values shown in the last line of table \ref{tab:PMwall} are related to the Sakiadis problem \cite{Sakiadis:1961:BLBa,Sakiadis:1961:BLBb} and were found by a few trial and miss attempts.
For this problem, the obtained skin friction coefficient is in good agreement with the values reported by other authors, e.g. $-0.44375$ Sakiadis \cite{Sakiadis:1961:BLBa}, $-0.4438$ by Ishak et al. \cite{Ishak:2007:BLM}, $-0.44374733$ by Cortell \cite{Cortell:2010:NCB} or $-0.443806$ by Fazio \cite{Fazio:2015:ITM}.

\begin{figure}[!hbt]
	\centering
\psfrag{es}[][]{$\eta^*$} 
\psfrag{e}[][]{$\eta$} 
\psfrag{df}[l][]{$\frac{df}{d\eta}$} 
\psfrag{ddf}[l][]{$\frac{d^2f}{d\eta^2}$} 
\psfrag{df*}[l][]{$\frac{df^*}{d\eta^*}$} 
\psfrag{ddf*}[l][]{$\frac{d^2f^*}{d{\eta^*}^2}$} 
\includegraphics[width=\www,height=6cm]{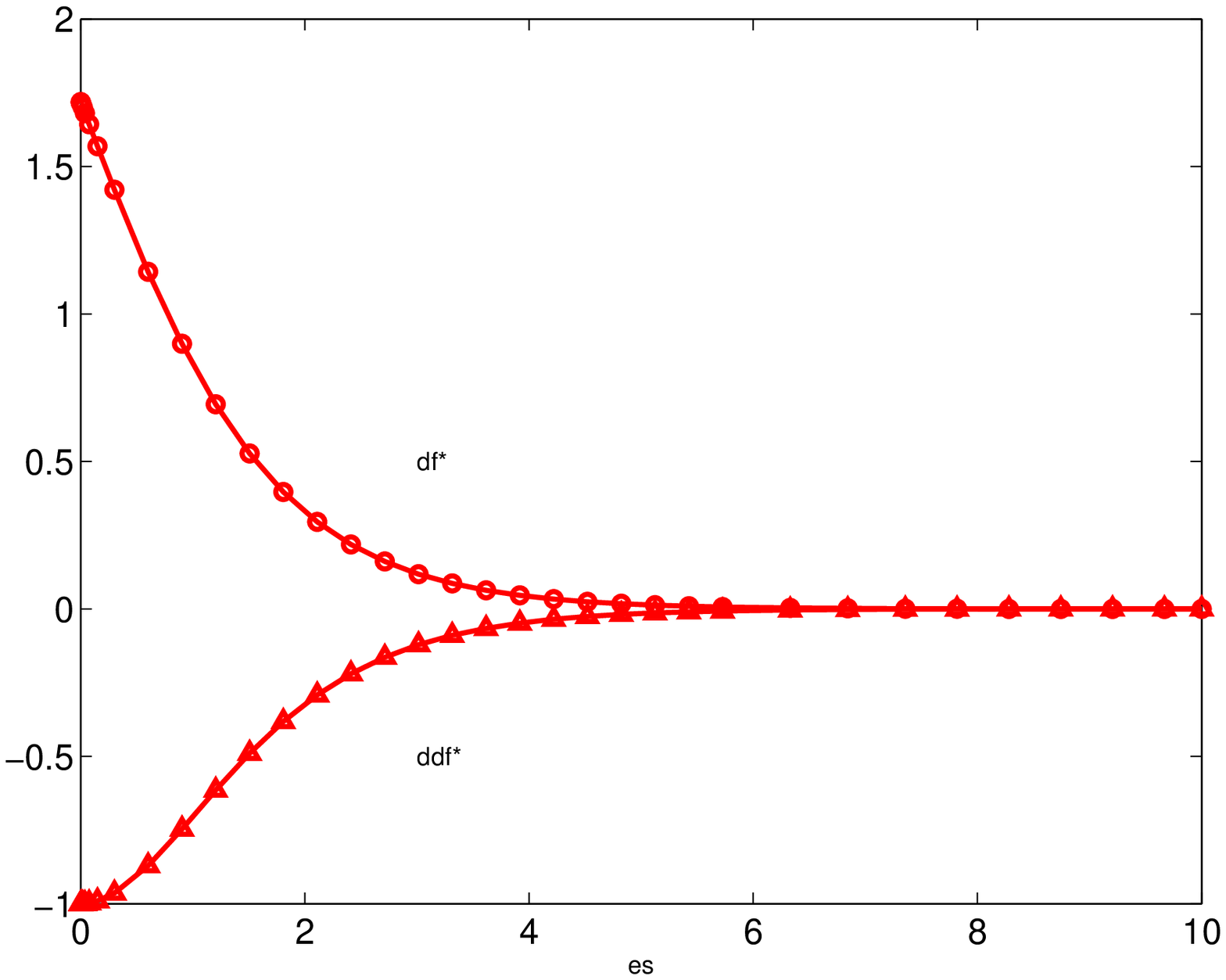} \\
\includegraphics[width=\www,height=6cm]{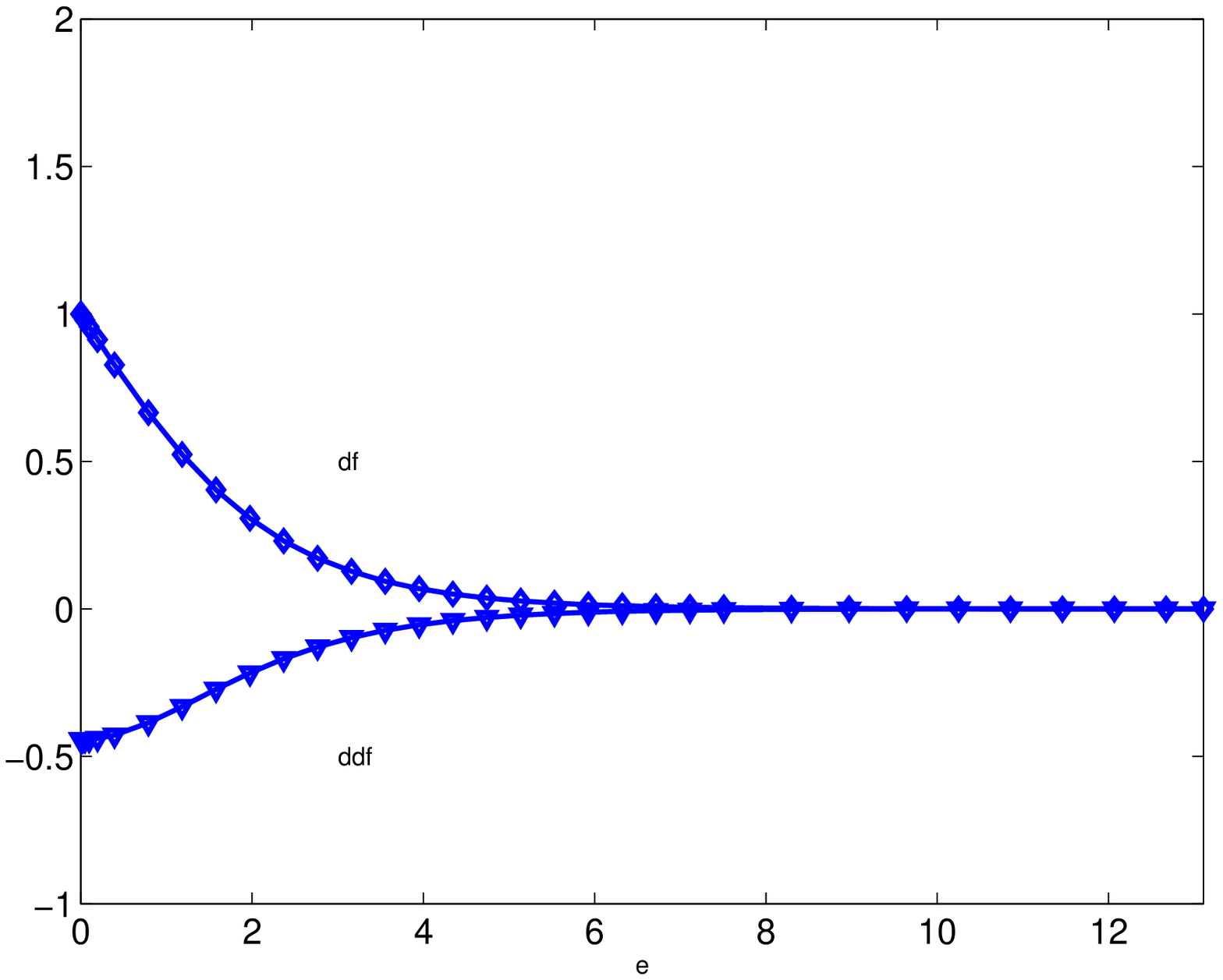}
\caption{Numerical results of the non-ITM. Top frame: solution of the IVP; bottom frame: solution of the Sakiadis problem found after rescaling.}
	\label{fig:Sakiadis}
\end{figure}
Figure \ref{fig:Sakiadis} shows the solution of the Sakiadis problem, describing the behaviour of a boundary layer flow due to a moving flat surface immersed in an otherwise quiescent fluid, corresponding to $P=1$.
Actually, this is a case of practical interest if we are considering the plate as an idealization of an airplane wing.
Let us notice here that by rescaling we get $\eta^*_{\infty} < \eta_{\infty}$. 

\begin{figure}[!hbt]
	\centering
\psfrag{r}[][]{$P$} 
\psfrag{df-}[][r]{$\ \ \ \ \frac{df^*}{d\eta^*}(0)=-1$} 
\psfrag{df+}[l][l]{\ \ \ $\frac{df^*}{d\eta^*}(0)=1$} 
\psfrag{d2f}[][]{$\frac{d^2f}{d\eta^2}(0)$} 
\includegraphics[width=\www]{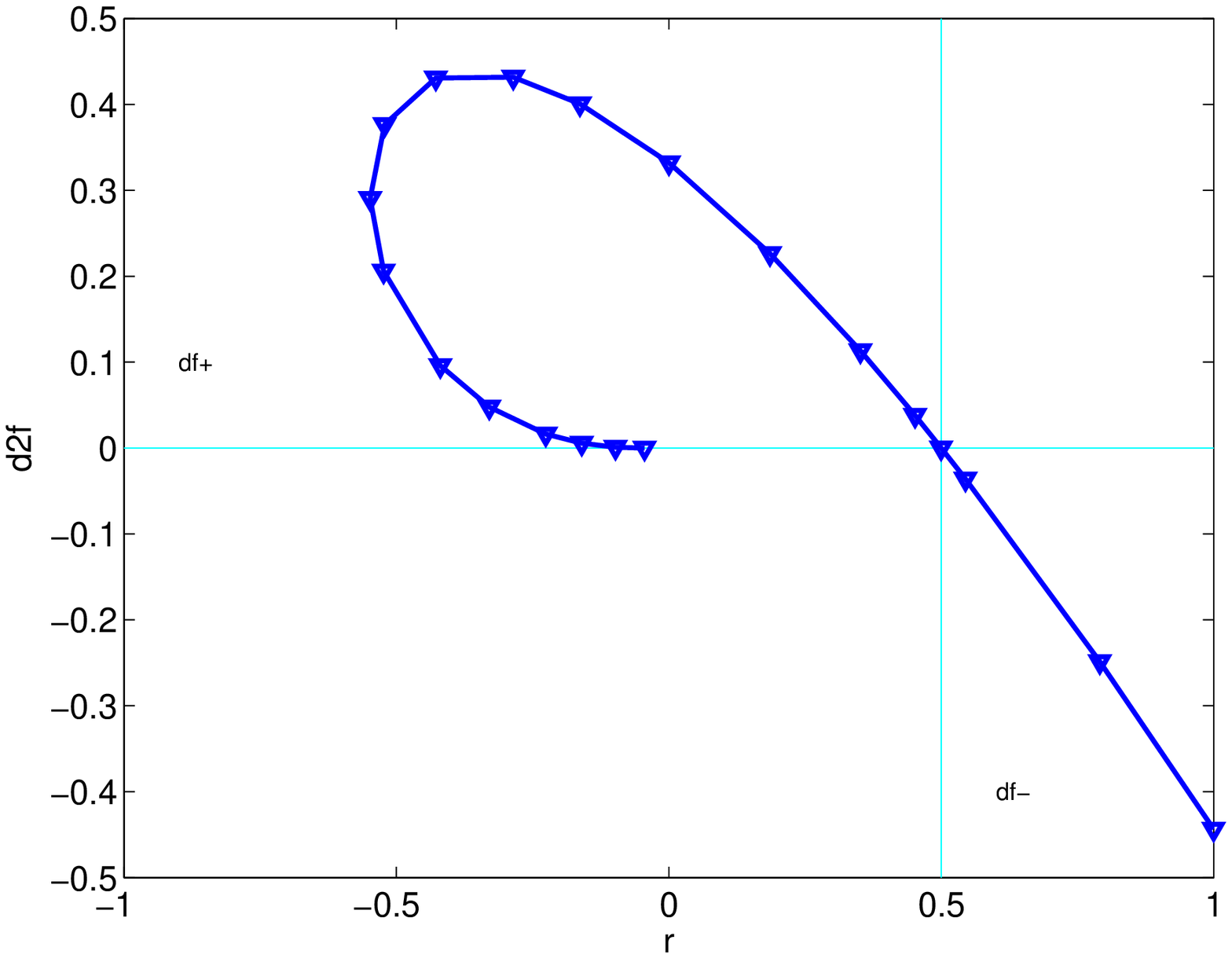}
\caption{Blasius problem with moving wall boundary conditions. Non-ITM: $\frac{d^2f}{d\eta^2}(0)$ versus $P$.}
	\label{fig:P-Mwall}
\end{figure}
In figure \ref{fig:P-Mwall} we plot $\frac{d^2f}{d\eta^2}(0)$ versus $P$.
From this figure we realize that our problem has an unique solution for $P\ge 0$, while dual solutions exists for $P<0$.
The values of the skin friction coefficient are positive for $P<1/2$ and they become negative when $P>1/2$.
From a physical viewpoint, this means that a drag force is exerted by the flow on the plate in the first case, while in the second the force is of opposite type.
Figure \ref{fig:P-Mwall} also shows that the solutions exist until a critical, negative, value of the parameter, say $P_c$, is reached.
From the data in table \ref{tab:PMwall} we get $P_c \approx -0.548447$.   
This value is in good agreement with the value $P_c = -0.5483$ computed by Ishak et al. \cite{Ishak:2007:BLM} using an iterative method: the second order Keller's Box finite difference method, see Keller \cite{Keller74,Keller}.
The boundary layer separate from the surface at $P=P_c$, and, therefore, for $P < P_c$ the Navier-Stoker equations has to be solved because the hypotheses of boundary layer theory felt down.
We have separation for a positive value of the skin friction coefficient and not at the point where this coefficient vanish as in the classical boundary layer theory. 
The zero value of the skin friction coefficient when $P=1/2$ corresponds to equal velocity of the plate and the free stream and does not mark separation.

\section{Surface gasification boundary condition}
In this section we consider a surface gasification flow that, with and without subsequent gas-phase flame-sheet burning, was first formulated and analysed by Emmons \cite{Emmons:1956:FCL}.
For such a flow we have to consider the variant of the celebrated Blasius problem
\begin{align}\label{eq:SGas}
& {\displaystyle \frac{d^3 f}{d \eta^3}} + f
{\displaystyle \frac{d^{2}f}{d\eta^2}}
= 0 \nonumber \\[-1ex]
&\\[-1ex]
& f(0) = -P \frac{d^2f}{d\eta^2} (0) \ , \qquad
\frac{df}{d\eta} (0) = 0 \ , \qquad
{\displaystyle \frac{df}{d\eta}}(\eta) \rightarrow 1 \quad \mbox{as}
\quad \eta \rightarrow \infty \ , \nonumber
\end{align}
where $ P $ is the classical Spalding heat transfer number \cite{Spalding:1961:MTT}. 
This transfer number for slow vaporization belongs to the interval $[0, 0.1]$ and varies from $P = O(1)$ to $P \approx 20$ for strong burning.
This problem has been studied recently by Lu and Law \cite{Lu:2014:ISB}.
These authors define an iterative method that has been shown to produce more accurate numerical results than the classical approximate solutions.

\subsection{The non-ITM} 
In the present case we consider the extended scaling group
\begin{equation}\label{eq:scalinv:gass1}
f^* = \lambda f \ , \qquad \eta^* = \lambda^{-1} \eta \ , \qquad 
P^* = \lambda^{-2} P \ .   
\end{equation}
Let us notice that the governing differential equation and the two boundary conditions at $\eta = 0$ in (\ref{eq:SGas}) are left invariant under (\ref{eq:scalinv:gass1}) and, on the contrary, the asymptotic boundary condition is not invariant.
By setting a value of $P^*$, we can integrate the Blasius governing differential equation in (\ref{eq:SGas}) in the star variables on $[0, \eta^*_\infty]$ with initial conditions
\begin{equation}\label{eq:ICs:gass}
f^*(0) = -P^* \ , \quad \frac{df^*}{d\eta^*}(0) = 0 \ , \quad \frac{d^2f^*}{d\eta^{*2}}(0) = 1 \ ,
\end{equation}
in order to compute $\frac{df^*}{d\eta^*}(\eta^*_\infty) \approx \frac{df^*}{d\eta^*}(\infty)$.
Here $\eta^*_\infty$ is a suitable truncated boundary.
The value of $ \lambda $ can be found by
\begin{equation}\label{eq:lambda4}
\lambda = \left[ \frac{d f^*}{d \eta^{*}}(\eta_\infty^*) \right]^{1/2} \ .   
\end{equation} 
After using (\ref{eq:lambda4}) to get the value of $\lambda$, we can apply the scaling invariance to obtain the missing initial conditions
\begin{equation}\label{eq:MICs:gas}
f(0) = \lambda^{-2} P^* \ , \quad \frac{d^2f}{d\eta^{2}}(0) = \lambda^{-3} \ .
\end{equation}

For the reader convenience, in table \ref{tab:P-Gas} we list sample numerical results.
\begin{table}[!htb]
\renewcommand\arraystretch{1.1}
	\centering
		\begin{tabular}{r@{.}lr@{.}lr@{.}lr@{.}lr@{.}l}
\hline \\[-3ex]
\multicolumn{2}{c}%
{$ P^* $}
& \multicolumn{2}{c}%
{${\displaystyle \frac{df^*}{d\eta^*}(\infty)}$}
& \multicolumn{2}{c}%
{${-f(0)}$}
& \multicolumn{2}{c}%
{${\displaystyle \frac{d^2f}{d\eta^2}(0)}$}
& \multicolumn{2}{c}%
{$ P $} \\[2ex]
\hline
0  &    & 1  & 655301 & 0 &        & 0 & 469553 & 0 &         \\
0  & 1  & 1  & 793644 & 0 & 074668 & 0 & 416289 & 0 & 179364  \\
0  & 25 & 2  & 025902 & 0 & 175643 & 0 & 346795 & 0 & 506476  \\
0  & 5  & 2  & 485809 & 0 & 317129 & 0 & 255152 & 1 & 242904  \\
0  & 75 & 3  & 048481 & 0 & 429556 & 0 & 187877 & 2 & 286361  \\
1  &    & 3  & 726397 & 0 & 518031 & 0 & 139016 & 3 & 726397  \\
1  & 25 & 4  & 528469 & 0 & 587401 & 0 & 103770 & 5 & 660586  \\
1  & 5  & 5  & 469166 & 0 & 641403 & 0 & 078184 & 8 & 203749  \\
1  & 75 & 6  & 548781 & 0 & 683845 & 0 & 059670 &11 & 460366  \\
2  &    & 7  & 779561 & 0 & 717055 & 0 & 046086 &15 & 559122  \\
2  & 2  & 8  & 863956 & 0 & 738939 & 0 & 037893 &19 & 500704  \\
\hline			
		\end{tabular}
	\caption{Surface gasification boundary condition: non-ITM results.}
	\label{tab:P-Gas}
\end{table}
The case $P^*=P=0$ is, again, the Blasius problem (\ref{eq:SGas}).
In this case our non-ITM reduces to the original method defined by T\"opfer \cite{Topfer:1912:BAB}.
The obtained skin friction coefficient is in good agreement with the values available in literature, see for instance the value $0.46599988361$ computed by Fazio \cite{Fazio:1992:BPF}. 
On the other hand, our value is different from the value $0.490$ obtained by a $3-2$ iteration solution of Lu and Law \cite{Lu:2014:ISB}.
For the numerical results reported here, depending on the behaviour of the numerical solution, we have used $\eta^*_\infty = 5$ or $\eta^*_\infty = 10$.

\begin{figure}[!hbt]
	\centering
\psfrag{e}[][]{$\eta^*$, $\eta$} 
\psfrag{f2star}[][]{$\frac{df^*}{d\eta^*}$} 
\psfrag{f2}[][]{$\frac{df}{d\eta}$} 
\psfrag{f3star}[][]{$\frac{d^2f^*}{d\eta^{*2}}$} 
\psfrag{f3}[][]{$\frac{d^2f}{d\eta^2}$} 
\includegraphics[width=\www]{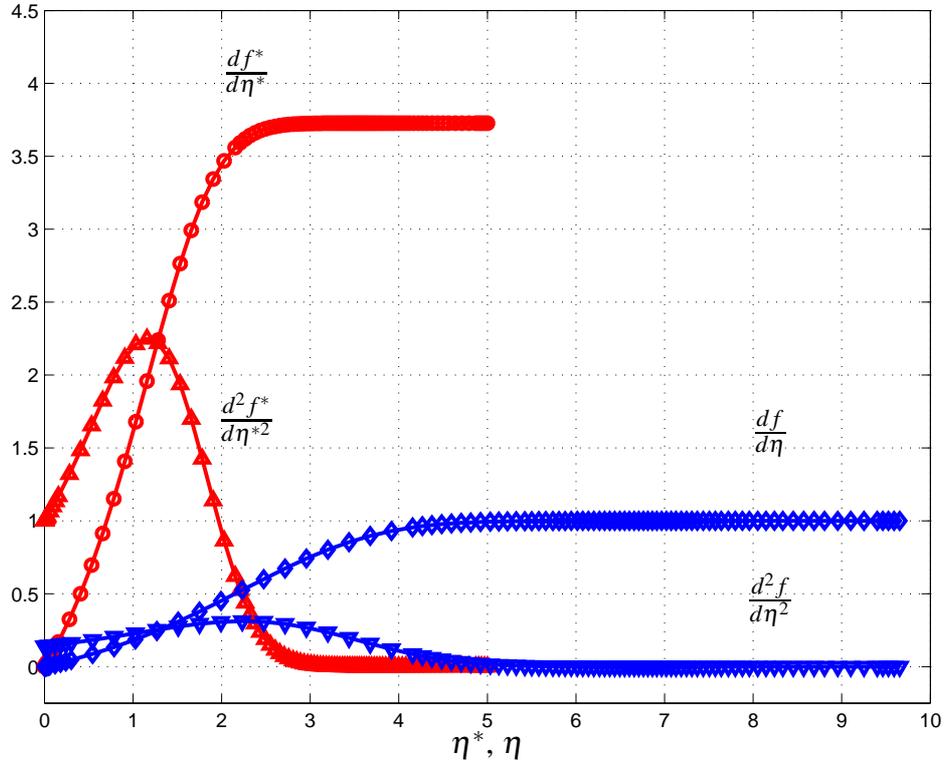}
\caption{Surface gasification boundary conditions with $P^* = 1$. Numerical solution by the non-ITM.}
	\label{fig:Blasius-gas}
\end{figure}

Figure \ref{fig:Blasius-gas} shows a sample numerical integration for $ P^* = 1 $ that is transformed under (\ref{eq:scalinv:gass1}) to $P \approx 3.726397$.
We notice that the solution of the Blasius problem with surface gasification boundary condition is computed by rescaling.
Moreover, by rescaling we get $\eta^*_{\infty} < \eta_{\infty}$. 

\begin{figure}[!hbt]
	\centering
\psfrag{P}[][]{$P$} 
\psfrag{fdf}[][]{$f(0)$, $\frac{d^2f}{d\eta^2}(0)$} 
\psfrag{f}[][]{$f(0)$} 
\psfrag{d2f}[][]{$\frac{d^2f}{d\eta^2}(0)$} 
\includegraphics[width=\www]{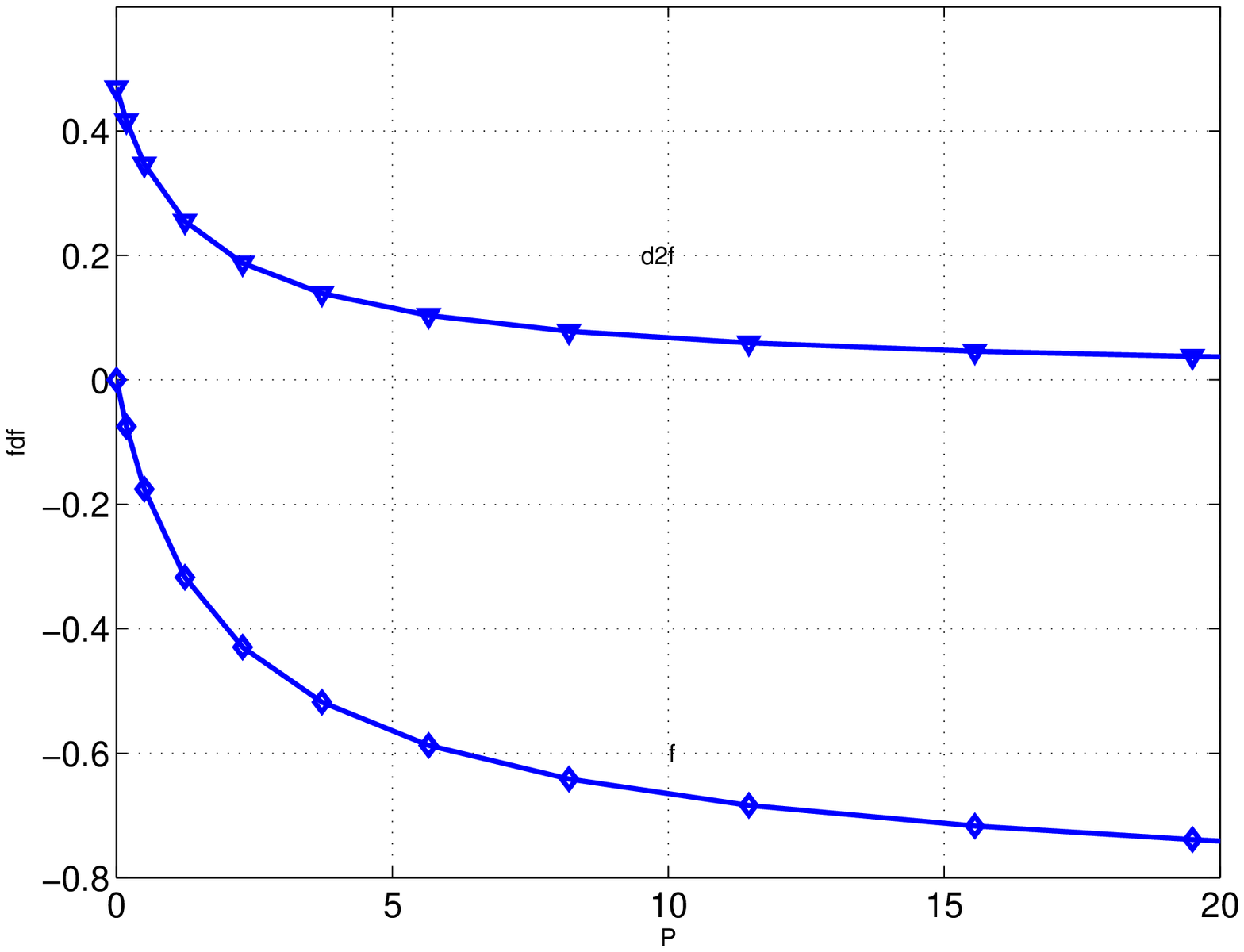}
\caption{Blasius problem with surface gasification boundary conditions. Non-ITM: behaviour of $f(0)$ and $\frac{d^2f}{d\eta^2}(0)$ with respect to $P$.}
	\label{fig:P-gas}
\end{figure}

In figure \ref{fig:P-gas} we plot $f(0)$ and $\frac{d^2f}{d\eta^2}(0)$ versus $P$.
We know, from the literature, that as $P$ goes to infinity then $f(0)$ goes to $-0.876$. 
Moreover, as it is easily seen, as $P$ goes to infinity then $\frac{d^2f}{d\eta^2}(0)$ goes to zero. 

\section{Concluding remarks}
The main contribution of this paper is the extension of the non-ITM proposed by T\"opfer \cite{Topfer:1912:BAB} for the numerical solution of the celebrated Blasius problem \cite{Blasius:1908:GFK} to classes of problems depending on a parameter.
By requiring the invariance of the involved parameter we are able to solve the given BVP non-iteratively but for a different value of the parameter.
This kind of extension was considered first by Na \cite{Na:1970:IVM}, see also NA \cite[Chapters 8-9]{Na:1979:CME}.
Here we defined a non-ITM for Blasius equation with moving wall or surface gasification. 

Finally, the reader should be advised that non-ITM cannot be applied to all problems of boundary layer theory.
In fact, sometimes we have to face problems that are not invariant with respect to all scaling groups.
As an example we can consider the Falkner-Skan model \cite{Falkner:1931:SAS}
\begin{align}\label{eq:abf}
& {\displaystyle \frac{d^{3}f}{d\eta^3}} + f 
{\displaystyle \frac{d^{2}f}{d\eta^2}} + P \; \left[1 - \left({\displaystyle
\frac{df}{d\eta}}\right)^2 \right] = 0 \ , \nonumber \\[-1.2ex]
& \\[-1.2ex]
& f(0) = {\displaystyle \frac{df}{d\eta}}(0) = 0 \ , \qquad
{\displaystyle \frac{df}{d\eta}}(\eta) \rightarrow 1 \quad \mbox{as}
\quad \eta \rightarrow \infty \ , \nonumber
\end{align}
where $f$ and $\eta$ are similarity variables and $P$ is a parameter related to the functional form od the fluid mainstream velocity.
If we test the invariance of the governing differential equation in (\ref{eq:abf}) under the extended scaling group
\begin{equation}\label{eq:sgext}
\eta^* = \lambda^{\alpha_1} \eta \ , \qquad f^* = \lambda^{\alpha_2} f \ , \qquad 
P^* = \lambda^{\alpha_3} P \ ,   
\end{equation}
where $\lambda$ is, again, the group parameter and ${\alpha_j}$, for $j = 1, 2, 3$, are constant to be determined, then we get three invariant conditions
\begin{equation}\label{eq:icond}
\alpha_2 -3 \alpha_1 = 2 (\alpha_2 - \alpha_1) = \alpha_3 = \alpha_3 + 2 (\alpha_2 - \alpha_1) \ .
\end{equation}
Now, it is a simple matter to show that the linear system defined by (\ref{eq:icond}) has the unique solution $\alpha_1 = \alpha_2 = \alpha_3 = 0$.
However, an iterative extension of our transformation method has been developed in \cite{Fazio:1994:FSE,Fazio:1996:NAN} and successfully applied to the Falkner-Skan model \cite{Fazio:1994:FSE,Fazio:2013:BPF}.

\vspace{1.5cm}

\noindent {\bf Acknowledgement.} {This research was 
partially supported by the University of Messina and by the GNCS of INDAM.}


\begin{thebibliography}{10}

\bibitem{Prandtl:1904:UFK}
L.~Prandtl.
\newblock {\"U}ber {F}l{\"u}ssigkeiten mit kleiner {R}eibung.
\newblock In {\em Proceedings Third Internernatinal Math. Congress}, pages
  484--494, 1904.
\newblock Engl. transl. in {NACA} Tech. Memo. 452.

\bibitem{Blasius:1908:GFK}
H.~Blasius.
\newblock Grenzschichten in {F}l\"{u}ssigkeiten mit kleiner {R}eibung.
\newblock {\em Z. Math. Phys.}, 56:1--37, 1908.

\bibitem{Schlichting:2000:BLT}
H.~Schlichting and K.~Gersten.
\newblock {\em Boundary Layer Theory}.
\newblock Springer, Berlin, 8th edition, 2000.

\bibitem{Weyl:1942:DES}
H.~Weyl.
\newblock On the differential equation of the simplest boundary-layer problems.
\newblock {\em Ann. Math.}, 43:381--407, 1942.

\bibitem{Topfer:1912:BAB}
K.~T{\"o}pfer.
\newblock Bemerkung zu dem {A}ufsatz von {H}. {B}lasius: {G}renzschichten in
  {F}l{\"u}ssigkeiten mit kleiner {R}eibung.
\newblock {\em Z. Math. Phys.}, 60:397--398, 1912.

\bibitem{Serrin:1970:ETS}
J.~Serrin.
\newblock Existence theorems for some compressible boundary layer problems.
\newblock In {\em Proc. of the Conference on Qualitative Theory of Nonlinear
  Differential and Integral Equations}, volume~5 of {\em SIAM studies in
  Applied Mathematics}, pages 35--42, 1970.

\bibitem{Meyer:1971:IMF}
R.~E. Meyer.
\newblock {\em Introduction to Mathematical Fluid Dynamics}.
\newblock Wiley, New York, 1971.

\bibitem{Hastings:2012:CMO}
S.~P. Hastings and J.~B. McLeod.
\newblock {\em Classical Methods in Ordinary Differential Equations With
  Applications to Boundary Value Problems}, volume 129 of {\em Graduate Studies
  in Mathematics}.
\newblock American Mathematical Society, Providence, 2012.

\bibitem{Rubel:1955:EET}
L.~A. Rubel.
\newblock An estimation of the error due to the truncated boundary in the
  numerical solution of the {B}lasius equation.
\newblock {\em Quart. Appl. Math.}, 13:203--206, 1955.

\bibitem{Fazio:2002:SFB}
R.~Fazio.
\newblock A survey on free boundary identification of the truncated boundary in
  numerical {BVP}s on infinite intervals.
\newblock {\em J. Comput. Appl. Math.}, 140:331--344, 2002.

\bibitem{Na:1970:IVM}
T.~Y. Na.
\newblock An initial value method for the solution of a class of nonlinear
  equations in fluid mechanics.
\newblock {\em J. Basic Engrg. Trans. ASME}, 92:503--509, 1970.

\bibitem{Na:1979:CME}
T.~Y. Na.
\newblock {\em Computational Methods in Engineering Boundary Value Problems}.
\newblock Academic Press, New York, 1979.

\bibitem{Gad-el-Hak:1999:FMM}
M.~Gad el~Hak.
\newblock The fluid mechanics of microdevices --- the {F}reeman scholar
  lecture.
\newblock {\em J. Fluids Eng.}, 121:5--33, 1999.

\bibitem{Martin:2001:BBL}
M.~J. {Martin} and I.~D. {Boyd}.
\newblock {Blasius boundary layer solution with slip flow conditions}.
\newblock In {\em Rarefied Gas Dynamics: 22nd International Symposium}, volume
  585 of {\em American Institute of Physics Conference Proceedings}, pages
  518--523, 2001, DOI: 10.1063/1.1407604.

\bibitem{Fazio:2009:NTM}
R.~Fazio.
\newblock Numerical transformation methods: {B}lasius problem and its variants.
\newblock {\em Appl. Math. Comput.}, 215:1513--1521, 2009.

\bibitem{Ishak:2007:BLM}
A.~Ishak, R.~Nazar, and I.~Pop.
\newblock Boundary layer on a moving wall with suction and injection.
\newblock {\em Chin. Phys. Lett.}, 24:2274--2276, 2007.

\bibitem{Emmons:1956:FCL}
H.~W. Emmons.
\newblock The film combustion of liquid fluid.
\newblock {\em ZAMM - J. Appl. Math. Mech.}, 36:60--71, 1956.

\bibitem{Lu:2014:ISB}
Z.~Lu and C.~K. Law.
\newblock An iterative solution of the {B}lasius flow with surface
  gasification.
\newblock {\em Int. J. Heat and Mass Transfer}, 69:223--229, 2014.

\bibitem{Sakiadis:1961:BLBa}
B.~C. Sakiadis.
\newblock Boundary-layer behaviour on continuous solid surfaces: I.
  {B}oundary-layer equations for two-dimensional and axisymmetric flow.
\newblock {\em AIChE J.}, 7:26--28, 1961.

\bibitem{Sakiadis:1961:BLBb}
B.~C. Sakiadis.
\newblock Boundary-layer behaviour on continuous solid surfaces: {II}. {T}he
  boundary layer on a continuous flat surface.
\newblock {\em AIChE J.}, 7:221--225, 1961.

\bibitem{Fazio:2015:ITM}
R.~Fazio.
\newblock The iterative transformation method for the {S}akiadis problem.
\newblock {\em Comput. \& Fluids}, 106:196--200, 2015.

\bibitem{Fazio:1992:BPF}
R.~Fazio.
\newblock The {Blasius} problem formulated as a free boundary value problem.
\newblock {\em Acta Mech.}, 95:1--7, 1992.

\bibitem{Boyd:1999:BFC}
J.~P. Boyd.
\newblock The {B}lasius function in the complex plane.
\newblock {\em Exp. Math.}, 8:381--394, 1999.

\bibitem{Cortell:2010:NCB}
R.~Cortell Bataller.
\newblock Numerical comparisons of {B}lasius and {S}akiadis flows.
\newblock {\em MATEMATIKA}, 26:187--196, 2010.

\bibitem{Keller74}
H.~B. Keller.
\newblock Accurate difference methods for nonlinear two-point boundary value
  problems.
\newblock {\em SIAM J. Numer. Anal.}, 11:305--320, 1974.

\bibitem{Keller}
H.~B. Keller.
\newblock {\em Numerical Methods for Two-point Boundary Value Problems}.
\newblock Dover Publications, New York, 2nd edition, 1992.

\bibitem{Spalding:1961:MTT}
D.~B. Spalding.
\newblock Mass transfer through laminar boundary layers - 1. {T}he velocity
  boundary layer.
\newblock {\em Int. J. Heat Mass Transfer}, 2:15--32, 1961.

\bibitem{Falkner:1931:SAS}
V.~M. Falkner and S.~W. Skan.
\newblock Some approximate solutions of the boundary layer equations.
\newblock {\em Philos. Mag.}, 12:865--896, 1931.

\bibitem{Fazio:1994:FSE}
R.~Fazio.
\newblock The {F}alkner-{S}kan equation: numerical solutions within group
  invariance theory.
\newblock {\em Calcolo}, 31:115--124, 1994.

\bibitem{Fazio:1996:NAN}
R.~Fazio.
\newblock A novel approach to the numerical solution of boundary value problems
  on infinite intervals.
\newblock {\em SIAM J. Numer. Anal.}, 33:1473--1483, 1996.

\bibitem{Fazio:2013:BPF}
R.~Fazio.
\newblock {B}lasius problem and {F}alkner-{S}kan model: {T}{\"o}pfer's
  algorithm and its extension.
\newblock {\em Comput. \& Fluids}, 73:202--209, 2013.

\end{thebibliography}
\end{document}